\documentclass[11pt,a4paper,leqno,twoside, headinclude]{amsart}

\usepackage[tracking=true]{microtype}
\DeclareMicrotypeSet*[tracking]{my}%
  { font = */*/*/sc/* }%
\SetTracking{ encoding = *, shape = sc }{ 45 }

\title[The Omnibus Conjecture---disproved]{The Omnibus Conjecture---disproved}
\def\titl{The Omnibus Conjecture---disproved}
\def\auth{Manuel Amann}
\date{November 3rd, 2020}


\author{\auth}

\usepackage[english]{babel}

\usepackage[colorlinks,pdftex, plainpages=false]{hyperref}
\hypersetup{pdftitle=\titl, pdfauthor=\auth, pdftoolbar=false,
plainpages=false, hyperindex=true, pdfdisplaydoctitle=true}

\hypersetup{colorlinks,%
citecolor=black,%
filecolor=black,%
linkcolor=black,%
urlcolor=black,%
pdftex}

\usepackage{color}

\usepackage{amsmath, amssymb, amscd, amsthm}
\usepackage{stmaryrd}
\usepackage{latexsym}
\usepackage[all]{xy}
\usepackage{pb-diagram}
\usepackage{rotating}
\usepackage{multicol}
\usepackage{lscape}
\usepackage{wasysym}
\usepackage{longtable}
\usepackage{enumerate}

\xyoption{all}

\newtheorem{theo}{Theorem}[section]
\newtheorem{main}{Theorem}

\newtheorem*{main*}{Theorem}

\newtheorem*{mainprop*}{Proposition}

\newtheorem{mainconj}{Conjecture}
\newenvironment{mainconjec}[1]{\begin{mainconj}[#1]\normalfont}{\end{mainconj}}

\newtheorem{prop}[theo]{Proposition}
\newtheorem{defi2}[theo]{Definition}
\newtheorem*{defi2*}{Definition}
\newenvironment{defi}{\begin{defi2}\normalfont}{\end{defi2}}
\newenvironment{defi*}{\begin{defi2*}\normalfont}{\end{defi2*}}

\newenvironment{defin*}[1]{\begin{defi2*}[#1]\normalfont}{\end{defi2*}}
\newtheorem*{rem2*}{Remark}
\newenvironment{rem*}{\begin{rem2*}\normalfont}{\hfill$\boxbox$\end{rem2*}}
\newtheorem{rem2}[theo]{Remark}
\newenvironment{rem}{\begin{rem2}\normalfont}{\hfill$\boxbox$\end{rem2}}

\newtheorem{lemma}[theo]{Lemma}
\newtheorem{cor}[theo]{Corollary}
\newtheorem*{cor*}{Corollary}

\newtheorem{conj}[theo]{Conjecture}
\newtheorem{ques}[section]{Question}
\newtheorem*{conj*}{Conjecture}
\newtheorem*{theo*}{Theorem}
\newtheorem*{ques*}{Question}
\newtheorem*{mi2}{Main Idea}

\newtheorem{ex2}[theo]{Example}
\newenvironment{ex}{\begin{ex2}\normalfont}{\hfill$\boxbox$\end{ex2}}

\newtheorem{exer2}[theo]{Exercise}

\newtheorem{alg2}[theo]{Algorithm}
\newenvironment{alg}{\begin{alg2}\normalfont}{\hfill$\boxbox$\end{alg2}}

\newtheorem{constr2}[theo]{Construction}
\newenvironment{constr}{\begin{constr2}\normalfont}{\hfill$\boxbox$\end{constr2}}

\newcommand{\cc}{{\mathbb{C}}}                                     
\newcommand{\nn}{{\mathbb{N}}}                                     
\newcommand{\qq}{{\mathbb{Q}}}                                     
\newcommand{\pp}{{\mathbf{P}}}                                     
\newcommand{\s}{{\mathbb{S}}}                                      
\newcommand{\zz}{{\mathbb{Z}}}                                     
\newcommand{\dif} {{\operatorname{d}}}                             
\newcommand{\In} {{\,\subseteq\,}}                                 
\newcommand{\im} {{\operatorname{im\,}}}                           
\newcommand{\Hom}{{\operatorname{Hom}}}                            
\newcommand{\APL}{{\operatorname{A_{PL}}}}                         
\newcommand{\Hur}{{\operatorname{Hur}}}                            
\newcommand{\ad}{{\operatorname{ad}}}                              
\newcommand{\co}{\colon\thinspace}                                 

\newcommand{\comment}[1]{}                                         
\newcommand{\xto}[1]{\xrightarrow{#1}}                             
\newcommand{\hto}[1]{\overset{#1}{\hookrightarrow}}                

\newcommand{\ack}{\noindent\textbf{Acknowledgements. }}            
\newcommand{\str}{\noindent\textbf{Structure of the article. }}    
\newcommand{\odd}{\textrm{odd}}                                    
\newcommand{\even}{\textrm{even}}                                  

\newenvironment{prf}{\begin{proof}[\textsc{Proof}]} {\end{proof}}     

\begin{document}

\maketitle \thispagestyle{empty}


\begin{abstract}
We provide various counter-examples to the long-standing so-called ``Omnibus Conjecture'' in Rational Homotopy Theory. That is, we show that a space with finite dimensional even-degree rational cohomology and finite dimensional spherical rational homology may indeed have infinite dimensional rational cohomology.

Moreover, we also discuss ``dual'' versions and special cases of the conjecture.
\end{abstract}


\section*{Introduction}

Rational Homotopy Theory can look upon a long success story. It allows for efficient computations of rational invariants of (nilpotent) topological spaces, and it reaches out to a multitude of applications in Riemannian and symplectic geometry. The theory has made huge progress since its foundations by Quillen and Sullivan. Nonetheless, it still leaves us with various prominent conjectures like the ``Halperin conjecture'', the ``toral rank conjecture'', or, the ``Hilali conjecture'' which resonate massively not only within the theory, but also within commutative algebra and geometry. In the seminal book \cite{FHT01} by F\'elix--Halperin--Thomas the ``Halperin conjecture'' and the ``toral rank conjecture'' (also due to Steve Halperin) are the first two open problems mentioned in a prominent collection of ``Seventeen Open Problems'' in the final chapter (see \cite[Chapter 39, p.~516]{FHT01}). Number 5 in this list is
\begin{mainconjec}{Omnibus Conjecture}\label{conj01}
Suppose $X$ is a simply-connected space with rational homology of finite type. If $H^\even (X; \qq)$ and the image of the Hurewicz homomorphism are both finite dimensional, does it follow that $H^*(X;\qq)$  is finite dimensional?
\end{mainconjec}
It owes its name and goes back to a failed proof attempt by Steve Halperin on an multi-stop train ride to Louvain-la-Neuve in an ``omnibus train'' in 1981.

The conjecture is trivially true for formal spaces (as in this case any generator of the cohomology algebra is spherical). Moreover, due to Baues' theorem (see \cite{Bau77}, cf.~\cite[Example, p.~328]{FHT01}), it is known for spaces with vanishing even-degree cohomology (as in this case the space rationally is a one-point union of spheres)---see Section \ref{sec04} for a more elaborate discussion of these results.

Recall that the image of the rational Hurewicz homomorphism,  $\im \operatorname{Hur}=\operatorname{Spher}_*(X)\In H_*(X;\qq)$, is referred to as \emph{spherical homology}, i.e.~homology represented by homotopy groups. We refer the reader to Section \ref{sec01} for a detailed discussion of spherical (co)homology. Recall that we call a space of \emph{finite type} if so is its rational homology, i.e.~if for each $i\geq 1$ it holds that $\dim H_*(X;\qq)<\infty$.

\bigskip

The goal of this article is to produce various counter-examples to the Omnibus Conjecture. They will be provided in Construction \ref{constr01} leading to
\begin{main}[Omnibus crashed]\label{theoA}
For every $n\geq 2$ there are infinitely many arbitrarily highly connected coformal topological spaces $X$ of finite type with $\operatorname{Spher}_*(X)$ dual to $H^\even(X;\qq)$ of dimension $n$, and satisfying that $\dim H^*(X;\qq)=\infty$.
\end{main}
It is a trivial observation that the theorem for arbitrary $n\geq 2$ follows from the respective observation for $n=2$ just by taking one-point unions with even-dimensional spheres of sufficiently high dimensions. Taking Cartesian products with even-dimensional spheres or complex projective spaces provides examples $X$ with $\dim H^*(X;\qq)=\infty$, and $n=\dim \operatorname{Spher}_*(X)< \dim H^\even(X;\qq)<\infty$. Taking one-point unions respectively products with spheres certainly again produces coformal spaces. Indeed, spheres are coformal spaces, and in the first case a Lie model for the one-point union in non-negative degrees results as the free product of the respective Lie models (see \cite[Example 2, p.~331]{FHT01}). As for products it is easy to observe that on the dual side of Sullivan models all differentials are quadratic, which is equivalent to coformality (see \cite[Example 7, p.~334]{FHT01}).

As a consequence, for the rest of this article we merely focus on $n=2$. Clearly, such constructions allow for a multitude of further examples. 
Just to point out, by referring to ``infinitely many'' examples we speak of infinitely many rationally inequivalent examples. As the proof shows, already the examples for $n=2$ can be constructed for several different choices of additional parameters yielding various infinite series---see Remark \ref{rem02}.

\bigskip

The main conceptual idea for constructing the counter-examples is to extend free graded Lie algebras (of two generators) by an infinite sequence of generators of even degrees with non-trivial differentials all mapping into the free Lie algebra---we tried to provide a concise description in Construction \ref{constr01}. For this recall that Rational Homotopy Theory may equally be described by either using Sullivan algebras or via the first approach due to Quillen in terms of differential graded Lie algebras. For us the latter bear the advantage that the free Lie model of a topological space $X$ is generated by the rational homology $H_*(X;\qq)$ of the space, and the homology of the Lie model actually corresponds to $\pi_*(X)\otimes \qq$ up to suspension. This gives us better control on homology than an equivalent but dual approach using Sullivan models. (Drawing upon universal coefficients we may work with rational homology instead of rational cohomology in (order to refute) Conjecture \ref{conj01}.)

\bigskip

In Section \ref{sec04} we discuss ``dual'' notions of the Omnibus Conjecture, of known results and of the counter-examples. More precisely, in view of the duality of free connected chain Lie algebras and Sullivan algebras we interchange respective restrictions on homology and homotopy. For example, we suggest and as well easily answer in the negative (see Example \ref{ex01}) the  following ``dual'' of Conjecture \ref{conj01}: ``Let $X$ be a simply-connected space with rational homology of finite type. If $\dim \pi_\odd(X)\otimes \qq<\infty$ and $\dim \operatorname{Spher}^*(X)<\infty$. Does then $\dim \pi_*(X)\otimes \qq<\infty$ hold true?''

\bigskip

We conclude the article by presenting a discussion of a specialised conjecture, which was independently communicated to us by Aniceto Murillo and Yves F\'elix, and which we named
\begin{conj}[elliptic Omnibus Conjecture]\label{conj02}
Every simply-connected topological space of finite type with finite dimensional rational homotopy and finite dimensional \emph{even-degree} rational cohomology is rationally elliptic.
\end{conj}
Recall for this that a simply-connected space $X$ is called \emph{rationally elliptic} if both $\dim \pi_*(X)\otimes \qq$ and $\dim H^*(X;\qq)$ are finite dimensional.
We shall see that this conjecture, this special case of the Omnibus conjecture, in some way is closely related to Anick's famous conjecture (see Section \ref{sec06}).

\bigskip

\str In Section \ref{sec01} we provide necessary preliminaries to the techniques used throughout the article. In particular, we provide a short introduction to the theory of differential graded Lie algebras, and we discuss different descriptions of spherical (co)homology. In Section \ref{sec03} we then prove Theorem \ref{theoA}. That is, we construct several counter-examples to the Omnibus Conjecture. Section \ref{sec04} is devoted to ``dual'' notions of the Omnibus Conjecture and to known results thereof. In Section \ref{sec06} we specialise to a discussion of the elliptic Omnibus Conjecture and of its interdependencies with Anick's conjecture.

\bigskip

\ack The author is grateful to Aniceto Murillo and Yves F\'elix for commenting on a previous version of the article and for communicating Conjecture \ref{conj02}.

The author was supported both by a Heisenberg grant and his research grant AM 342/4-1 of the German Research Foundation; he is moreover a member of the DFG Priority Programme 2026.


\section{Preliminaries}\label{sec01}

By no means is this section intended or able to provide an introduction to Rational Homotopy Theory. Let us merely review some classical concepts which are relevant to us. As a general reference for Rational Homotopy Theory we refer the reader to \cite{FHT01}. We shall provide several references to the relevant results and concepts in this textbook; however, drawing on a certain familiarity of the reader with the basic concepts of the theory.

Throughout this article $X$ will denote a simply-connected topological space of finite type. Up to spatial realisation (see Section \ref{subsec03}) we may equivalently encode and describe the rational type of $X$ via a (minimal) Sullivan model $m\co (\Lambda V,\dif)\xto{\simeq} \APL(X)$ respectively via a (minimal) free Lie model. This duality will be reflected upon and exploited in several places throughout the article and, in particular, in this section.

(Co)homology will always be taken with rational coefficients; all vector spaces, Lie algebras, etc.~are over the rationals. All algebras and spaces are of finite type.

\subsection{Homotopy groups of Sullivan algebras}\label{subsec04}

For an introduction to the theory of Sullivan algebras we draw on \cite[Chapter 12, p.~138]{FHT01}.
For such an algebra $(\Lambda V,\dif)$ we recall the decomposition $\dif=\dif_0+\dif_1+\dif_2+ \ldots$ uniquely determined by $\dif_i\co V\to \Lambda^{i}V$ raising word length in $V$ by $i$. The differential $\dif_1$ is the ``quadratic part'' of the differential encoding Whitehead products and hence the Lie bracket of the corresponding differential graded Lie algebra---cf.~\cite[Definition, p.~176]{FHT01}. Minimality hence is equivalent to $\dif_0=0$.

From \cite[Theorem 15.11, p.~208]{FHT01} we recall the main theorem of Rational Homotopy Theory stating that for a minimal Sullivan algebra $V^i$ for $i\geq 1$ is  isomorphic to the dual of $\pi_i(X)\otimes \qq$ provided that $X$ is simply-connected and of finite type. We extend this slightly in terms of notation by defining the homotopy groups of a Sullivan algebra $(\Lambda V,\dif)$ to be
\begin{align}\label{eqnhom}
\pi(\Lambda V,\dif):= H(V,\dif_0)
\end{align}
In case the model $(\Lambda V,\dif)$ is not minimal, the induced projection $\tilde \dif|_V\co V=\Lambda V/\Lambda^{\geq 2} V \to \Lambda V/\Lambda^{\geq 2} V=V$ still is a well-defined differential which identifies with $\dif_0$, the linear part of $\dif$, considered as a morphism $V\to V$. Up to duality we then obtain the graded isomorphism of vector spaces
\begin{align}\label{eqn08}
\pi(\Lambda V,\dif)\cong \pi_*(X)\otimes \qq
\end{align}
which clearly extends the case for minimal $(\Lambda V,\dif)$ (which satisfies  $\pi(\Lambda V,\dif=\dif_{\geq 1})=H(V,0)=V$). Indeed, for this we recall that up to isomorphism we may split any Sullivan algebra $(\Lambda V,\dif)$ as the product of a minimal one $(\Lambda V',\dif)$ and a contractible one $(\Lambda C,\dif)$ (see \cite[Proposition 14.9, p.~187]{FHT01}). As a consequence
\begin{align*}
\pi(\Lambda V,\dif)\cong \pi(\Lambda V',\dif) \oplus \pi(\Lambda C,\dif)\cong V'\cong  \pi_*(X)\otimes \qq
\end{align*}

\subsection{Lie algebras and word length}\label{subsec01}

The main technical arguments of this article lie in the field of differential graded Lie algebras. It is certainly impossible to review their theory at this place; for their structure and entire role in Rational Homotopy Theory we point to \cite[Part IV]{FHT01}, \cite{FHT15}, and \cite{Nei10}, which also serve us as a main reference for the latter theory. Let us merely recall some properties which will become important throughout the article thereby fixing notation.

From \cite[Chapter 21, p.~283, Chapter 22(e), p.~294]{FHT01} we recall
\begin{defi}A \emph{graded Lie algebra $L$} is a graded vector space $L=\{L_i\}_{i\in \zz}$ together with a \emph{Lie bracket}, i.e.~with a linear map of degree $0$, $L\times L\to L$, $x\times y\to [x,y]$ satisfying
\begin{itemize}
\item \emph{(anti-)commutativity}
\begin{align*}
[x,y]=-(-1)^{\deg x \deg y} [y,x]
\end{align*}
\item and the \emph{Jacobi-identity}
\begin{align*}
[x,[y,z]]=[[x,y],z]+(-1)^{\deg x \deg y} [y,[x,z]]
\end{align*}
\end{itemize}
A \emph{differential graded Lie algebra} is a graded Lie algebra together with a differential $\delta$ which is a derivation, i.e.~which satisfies
\begin{align*}
\delta([x,y])=[\delta x,y]+(-1)^{\deg x} [x,\delta y]
\end{align*}
We call a differential graded Lie algebra a \emph{connected chain Lie algebra} if $L=\{L_i\}_{i\in \nn}$, i.e~if it is positively graded.

A \emph{free differential graded Lie algebra} $\mathbb{L}_V$ over the graded vector space $V$ is the Lie subalgebra of the tensor algebra $TV$ (together with the commutator bracket) generated by $V$. A free connected chain Lie algebra $(\mathbb{L}_V, \delta)$ is minimal if the \emph{linear
part} $\delta_V\co V\to V$ of the differential $\delta$ (again defined by $\im (\delta-\delta_V)\In [\mathbb{L}_V,\mathbb{L}_V]$) is zero.
\end{defi}
We remark that free connected chain algebras play a role dual to Sullivan algebras; in particular, we can always find such minimal free models for connected chain Lie algebras and simply-connected finite type topological spaces uniquely determined up to isomorphism (see \cite[Chapter 24, p.~322]{FHT01}).

For a free graded Lie algebra $\mathbb{L}_{\langle x_1, \ldots, x_n\rangle}$ (over the vector space with the fixed basis $(x_i)_{1\leq i\leq n}$) we may define the word length of a monomial $x\in \mathbb{L}_{\langle x_1, \ldots, x_n\rangle}$ as the overall number of $x_i$ (with multiplicities) needed to write it as iterated brackets in the $x_i$ needed to express it. (By convention, the element $0$ has any word length; clearly, an element $x_i$ has word length $1$, a term $[x_i,x_j]$ has word length $2$.) We speak of the word length of an arbitrary element, once the element can be written as a sum of monomials of identical wordlengths. Note that this number is unaffected by both (anti-)commutativity and the Jacobi identity whence it is a well-defined invariant on the whole Lie algebra.
We use the notation
\begin{align*}
(\mathbb{L}_{\langle x_1, \ldots, x_n\rangle})^{\operatorname{wl}=k}
\end{align*}
in order to refer to the vector subspace of all elements of word length equal to $k$. (Analogously, ``$\leq k$'' is used for the vector subspace of those elements of word length at most $k$.)

In the same way we define the word length of a monomial $x$ (again well-defined) \emph{in a certain element $x_i$} as the number of $x_i$ appearing in a bracket representation of $x$. We denote the vector subspace of all elements of word length in $x_i$ equal to $k$ by
\begin{align*}
\mathbb{L}_{\langle x_1, \ldots, x_n\rangle}^{\operatorname{wl_{x_i}}=k}
\end{align*}

We recall some more standard terminology. By
\begin{align*}
\ad^j(x)(y)=\overbrace{[x,[x\ldots,[x}^j,y]\ldots ]
\end{align*}
we denote the \emph{adjoint map} consisting of $j$-fold iterative bracketing with $x$.

By $sL=L^{(+1)}$ we denote \emph{suspension}, that is a $(+1)$-degree shift $(sL)_i=L_{i-1}$. By $s^{-1}$ we correspondinly denote \emph{desuspension}.

\subsection{Spatial realisation of a differential graded Lie algebra}\label{subsec03}

It is always possible to realise a finite type connected chain Lie algebra $(L,\delta)$ as a simply-connected CW-complex. We illustrate this in two steps.

For this we first use the functor $C^*$ from differential graded Lie algebras to commutative differential graded algebras given by
\begin{align}\label{eqn05}
C^*(\mathbb{L}_V,\delta)=\Hom(C_*(\mathbb{L}_V,\delta),\qq)
\end{align}
where $C_*(\cdot)$ is the Cartan--Eilenberg--Chevalley construction (see \cite[Chapter 23(a), p.~313, Chapter 22(b), p.~301]{FHT01}).

From \cite[Lemma 23.1, Proposition 23.2, p.~314]{FHT01} we see that the connected chain Lie algebra $(L,\delta)$ dualises to a (usually not minimal) Sullivan algebra $(\Lambda V,\dif)$ (satisfying that $sL\cong V$). We apply this dualisation to $(L,\delta)$ and to a free connected cochain Lie algebra $(\mathbb{L}_W,\partial)$ and obtain graded vector space isomorphisms
\begin{align}\label{eqn04}
 sH(W,\delta_W)\oplus \qq\cong H(\Lambda V,\dif) \qquad \textrm{and} \qquad 
 H(L,\delta)\cong s^{-1}\pi(\Lambda V,\dif)
\end{align}
where $\delta_W$ is the induced linear parts of the differentials. (The first property up to duality is Isomorphism \cite[(24.3), p.~326]{FHT01} respectively \cite[Proposition 22.8, p.~306]{FHT01}, the second one is from \cite[Proposition 23.3, p.~318]{FHT01} in view of Definition \eqref{eqnhom}. For this, in the latter case, up to duality, we identify $H(L,\delta)$ with the homotopy Lie algebra and even obtain an isomorphism of graded Lie algebras.) We remark that these two properties are perfectly dual to one another; they merely express that up to degree shift taking homotopy of the Lie model corresponds to taking homology of the Sullivan model and vice versa. We just stick to established notation which conceals this slightly.

As a next step, we apply spatial realisation as a CW-complex for Sullivan algebras. For this see \cite[Chapter 17, p.~237]{FHT01}, wherein spatial realisation is attained by successive \emph{Sullivan realisation} $\langle \cdot \rangle$ from commutative cochain algebras to simplical sets and \emph{Milnor realisation} $|\cdot|$ from simplicial sets to CW-complexes. We denote this combined spatial realisation by the composition  $|\langle \cdot \rangle|$. We then apply it to $C^*(\mathbb{L}_{W},\delta)$ (see \eqref{eqn05}) and obtain the CW-complex
\begin{align*}
X=|\langle C^*(\mathbb{L}_{W}, \delta)\rangle|
\end{align*}
Due to \cite[Theorem 17.10, p.~250]{FHT01} the Sullivan algebra $C^*(\mathbb{L}_{W},\delta)$ is a model of $X$, and, consequently, the initial Lie algebra $(\mathbb{L}_{W}, \delta)$ becomes a free Lie model of $X$. Consequently, adapting Properties \eqref{eqn04} we obtain
\begin{align}\label{eqn06}
 H_{i-1}(W,\delta_W)\cong H_{i}(X) \qquad \textrm{and} \qquad 
 H_{i-1}(L,\delta)\cong \pi_{i}(X)\otimes \qq
\end{align}
for $i\geq 2$.

\subsection{Spherical (Co)homology}\label{subsec02} Let $(\Lambda V,\dif)$ denote a minimal Sullivan model of the simply-connected topological space $X$ of finite type. The dual, $\Hur^*$, of the rationalised Hurewicz homomorphism $\Hur\co \pi_*(X)\otimes \qq \to H_*(X;\qq)$ fits into the commutative square
\begin{align*}
\xymatrix{
H^{>0}(\Lambda V,\dif)  \ar[rr]^<<<<<<<<<<<{H(m)}_\cong
\ar[d]^\zeta &&H^*(X;\qq)\ar[d]^{\Hur^*}\\
V \ar[rr]^<<<<<<<<<<<<{\nu}_<<<<<<<<<<<<\cong
&& \Hom(\pi_*(X),\qq)
}
\end{align*}
where $\zeta$ is the map $H^{>0}(\Lambda V,\dif)\to V$ induced by taking the quotient with $\Lambda^{\geq 2} V$, and $\nu$ establishes the duality of homotopy groups
(see \cite[p.~173]{FHT01}). It follows that \emph{spherical cohomology}, i.e.~$\operatorname{Spher}^*(X):=\im \Hur^*$, is identified,
\begin{align*}
\operatorname{Spher}^*(X)\cong \ker (\bar \dif|_V)
\end{align*}
with the kernel of the map
\begin{align}\label{eqn01}
\bar \dif|_V\co V\to \Lambda V/\dif(\Lambda^{\geq 2} V)
\end{align}
induced by $\dif$ (using the minimality of the model, i.e.~the fact that non-trivial closed forms in $V$ are not exact). (We point out that the quotient $\Lambda V/\dif(\Lambda^{\geq 2} V)$ is a quotient of \emph{vector spaces}, i.e.~we do \emph{not} divide out the ideal generated by $\dif(\Lambda^{\geq 2} V)$.)

(Up to a suitable isomorphism of the minimal model, i.e.~potentially for another choice of minimal model within its isomorphism class, called \emph{normal form}, we may identify this with $\ker \dif|_V$.)

In case the model $(\Lambda V,\dif)$ is not minimal we observed in Section \ref{subsec04} and Isomorphism \eqref{eqn08} that
\begin{align*}
H(V,\dif_0)\cong \Hom(\pi_*(X), \qq)
\end{align*}
again inducing and making commutative the analog of the above diagramme, namely
\begin{align*}
\xymatrix{
H^{>0}(\Lambda V,\dif)  \ar[rr]^<<<<<<<<<<<{H(m)}_\cong
\ar[d]^\zeta &&H^*(X;\qq)\ar[d]^{\Hur^*}\\
H(V,\dif_0) \ar[rr]^<<<<<<<<<<<<{\nu}_<<<<<<<<<<<<\cong
&& \Hom(\pi_*(X),\qq)
}
\end{align*}

In this description spherical cohomology becomes
\begin{align*}
\operatorname{Spher}^*(X)\cong\ker (\bar \dif|_V)/\im \dif_0 \In V
\end{align*}
using the map $\bar \dif$ from \eqref{eqn01} induced by $\dif$ via projection as well as the linear part $\dif_0$ of $\dif$.

\bigskip

Spherical homology in the free Lie model $(\mathbb{L}_V,\delta)$ of $X$ is obtained as follows---note the perfectly dual descriptions from \eqref{eqn04}. From \cite[Proposition 24.4, p.~326]{FHT01} we cite the commutativity of the corresponding dual diagramme
\begin{align*}
\xymatrix{
sH(\mathbb{L}_V,\delta)  \ar[rr]^<<<<<<<<<<<{}_\cong
\ar[d]^{sH(\mu)} && \pi_*(X)\otimes \qq\ar[d]^{\Hur}\\
sH(V,\delta_V) \ar[rr]
^<<<<<<<<<<<{}_<<<<<<<<<<<\cong
&& H_{>0}(X;\qq)
}
\end{align*}
The morphism $sH(\mu)$ is the suspension of the morphism induced in homology by the surjective morphism $\mu\co (\mathbb{L}_V,\delta)\to (V,\delta_V)$ obtained by taking the quotient by the ideal $[\mathbb{L}_V,\mathbb{L}_V]$ (preserved by $\delta$) in $\mathbb{L}_V$ of all elements of bracket lengths at least $2$. Since $\mathbb{L}_V=V\oplus [\mathbb{L}_V,\mathbb{L}_V]$ this then maps to $sH(V,\delta_V)$  with $\delta_V$ the differential induced on $V$ by $\delta$.

Hence \emph{spherical homology} $\operatorname{Spher}_*(X)$, i.e.~the image of the rational Hurewicz homomorphism, $\im \Hur$, on the level of free Lie models identifies with the following vector subspace, namely with the suspension
\begin{align*}
\operatorname{Spher}_*(X)\cong s\big(\ker \bar \delta|_V /\im \delta_V  \big) \In s V
\end{align*}
where
\begin{align*}
\bar\delta|_V\co V \to \mathbb{L}_V/\delta([\mathbb{L}_V,\mathbb{L}_V])
\end{align*}
is the projection of $\delta$ to the \emph{vector space} quotient. (As above, we do \emph{not} consider the ideal generated by $\delta([\mathbb{L}_V,\mathbb{L}_V])$.)
Again, for a minimal free Lie model we get that $sH(V,\delta_V)\cong sV$, and spherical homology given by
\begin{align*}
\operatorname{Spher}_*(X)\cong s\big(\ker \bar \delta|_V \big) \In s V
\end{align*}
This yields a description perfectly dual to all we deduced for Sullivan algebras.

In both cases we fix the following observation:
\begin{rem}\label{rem01}
Let either $(\Lambda V,\dif)$ be a minimal Sullivan algebra, or let $(\mathbb{L}_V,\delta)$ be a minimal free differential graded Lie algebra. Then an element $0\neq v\in V$ (in each respective case) constitutes a non-trivial spherical (co)homology class if and only if $\dif v$ respectively $\delta v$ lies in the vector subspace $\dif(\Lambda^{\geq 2} V)$ in $(\Lambda V,\dif)$, respectively in the vector subspace $\delta([\mathbb{L}_V,\mathbb{L}_V])$ in $(\mathbb{L}_V,\delta)$.
\end{rem}


\section{Proof of Theorem \ref{theoA}}\label{sec03}

In this section we construct counter-examples to the Omnibus Conjecture. As noted in the introduction due to one-point unions with even-dimensional spheres, it suffices to construct infinitely many arbitrarily highly connected finite type coformal topological spaces $X$ each satisfying $\dim H^*(X)=\infty$, and $\dim \operatorname{Spher}_*(X)=\dim H^\even(X;\qq)=2$.

We shall achieve this by constructing a minimal free differential graded Lie algebra with the corresponding properties, i.e.~a Lie algebra which then can be realised as a topological space $X$  having the properties above (see Section \ref{subsec03}).

The starting point of our constructions will be a free differential graded Lie algebra $(\mathbb{L}_{\langle a,b\rangle},0)$ on two generators $a$, $b$ of \emph{odd} degree $\deg a,\deg b\geq 1$. Indeed, the final construction will only rely on both degrees being odd. That is, here already by changing degrees we can achieve that this obviously leads to infinitely many rationally inequivalent examples. As we shall see in the following many more parameters for the examples can be chosen---see Remark \ref{rem02}.

Based upon this free Lie algebra $(\mathbb{L}_{\langle a,b\rangle},0)$ (i.e.~$\delta a=\delta b=0$) in each example we shall add an infinite sequence of new generators $(x_i)_{i\geq 1}$ for which we carefully need to specify their respective differentials $\delta x_i$. Any such construction will then be realised as a topological space
\begin{align*}
X=|\langle C^*(\mathbb{L}_{\langle a, b, x_i\rangle_{i\geq 1}}, \delta)\rangle|
\end{align*}
given as the spatial realisation $|\langle \cdot \rangle|$ (see Section \ref{subsec03}) of the Sullivan algebra dual to the differential graded Lie algebra---the latter then clearly becomes a Lie model of $X$.

Hence the main/remaining challenge in the construction of $X$ and, consequently, in the proof of Theorem \ref{theoA} lies in constructing the differentials $\delta x_i$ (and picking the degrees $\deg x_i$) for $i\geq 1$.

\bigskip

For this we want to explicitly construct a graded Hall basis for the vector subspace of the free graded Lie algebra $\mathbb{L}_{\langle a,b\rangle}$ of all those brackets with word-length exactly $2$ in $b$ (see Section \ref{subsec01}). Again denote this vector space by $\mathbb{L}_{\langle a,b\rangle}^{\operatorname{wl}_b=2}$. We derive an iterative algorithm for doing so from \cite[Examples 8.7.4 and 8.7.5]{Nei10}, both drawing upon \cite[Proposition 8.7.3, p.~272]{Nei10} wherein it is stated that the kernel, $K_{i+1}$, in the short exact sequence \eqref{eqnseq} is again a free algebra, and which establishes a formula for the Euler--Poincar\'e series which then allows to show that constructed elements of generating sets are necessarily linearly independent.
\begin{alg}\label{alg01}
Starting with $i=1$, the free connected cochain Lie algebra $\mathbb{L}_{Z_i}$ over the graded vector space $Z_i$, and iterating over successive $i\geq 1$ let
\begin{align}\label{eqnseq}
0\to K_{i+1}\to \mathbb{L}_{Z_i}\to \mathbb{L}_{\langle x_i\rangle}\to 0
\end{align}
be an exact sequence of free graded Lie algebras with homogeneous $0\neq x_i\in Z_i$ induced by projecting away a homogeneous complement $Z_i\to \langle x_i\rangle$. Then $K_{i+1}=\mathbb{L}_{Z_{i+1}}$ is a free Lie algebra over the new vector space $Z_{i+1}$ (defined by this splitting). More precisely, this vector space can be determined as follows: choose a homogeneous direct sum decomposition $Z_i=\langle x_i\rangle \oplus W_i$. We have to make the following distinction depending on the parity of $x_i$.
\begin{itemize}
\item
If $\deg x_i$ is odd, then
\begin{align*}
Z_{i+1}=W_i\oplus [x_i,W_i] \oplus [x_i,x_i]
\end{align*}
\item
If $\deg x_i$ is even (implying $[x_i,x_i]=0$), then
\begin{align*}
Z_{i+1}&=\bigoplus_{j\geq 0} \ad^j(x_i)(W_i)\\
\\&= W_i \oplus [x_i,W_i] \oplus [x_i,[x_i,W_i]] \oplus [x_i,[x_i,[x_i,W_i]]] \oplus \ldots
\end{align*}
\end{itemize}
We continue with splitting off a new homogeneous element $x_{i+1}$ from $Z_{i+1}$, i.e.~by passing from $i\to i+1$. The split-off elements $(x_i)_{i\geq 1}$ yield a basis. This allows to construct a basis of the free Lie algebra $\mathbb{L}_{Z_1}$ with $\dim Z_1<\infty$ up to a certain word length in finitely many steps.
\end{alg}

As a trivial application of this algorithm or actually an observation underlying it we recall
\begin{lemma}\label{lemma01}
As graded vector spaces we obtain that
\begin{align*}
\dim (\mathbb{L}_{\langle a,b\rangle}^{\operatorname{wl}_b=0})=\langle a, [a,a]\rangle
\end{align*}
\end{lemma}
\begin{prf}
We do not apply the algorithm, but simply observe that via the Jacobi identity one may produce a generating set of $\dim (\mathbb{L}_{\langle a,b\rangle}^{\operatorname{wl}_b=0})$ via successive brackets of the form $[a,[a,[\ldots,a]\ldots ]$. That is any element of a certain word length is a multiple of this iterated element of the same word length. Since $\deg a$ is odd, another trivial application of the Jacobi identity yields that already a triple-product of this form vanishes.
\end{prf}

We now apply Algorithm \ref{alg01} in order to recursively provide a basis of $\mathbb{L}_{\langle a,b\rangle}^{\operatorname{wl}_b=2}$. In particular, it is our goal to show that $\dim (\mathbb{L}_{\langle a,b\rangle}^{\operatorname{wl}_b=2})^i$ is monotonously growing in $i$. (We shall need slightly more information though.)

In Table \ref{table01} the notation $[i,j]$ is shorthand notation for the element
\begin{align*}
[i,j]:=\big[[\overbrace{a,[a,[\ldots,[a,b}^{i}]\ldots], [\overbrace{a,[a,[\ldots,[a,b}^{j}]\ldots]\big]
\end{align*}
\begin{prop}\label{prop01}
A homogeneous basis of $\mathbb{L}_{\langle a,b\rangle}^{\operatorname{wl}_b=2}$ is provided in the respective degrees by Table \ref{table01}.

As a general formula we deduce that a basis of word length $n$ is given by
\begin{align*}
\begin{cases}
\{[1,n-1], [2,n-2], \ldots, [n/2-1,n/2+1]\} & \textrm{for } n\equiv 0 \mod 4\\
\{[1,n-1], [2,n-2], \ldots, [(n-1)/2,(n+1)/2]\} & \textrm{for } n\equiv 1,3 \mod 4\\
\{[1,n-1], [2,n-2], \ldots, [n/2,n/2]\}& \textrm{for } n\equiv 2 \mod 4
\end{cases}
\end{align*}
This yields its dimension in word length $n\geq 1$ as
\begin{align*}
\dim (\mathbb{L}_{\langle a,b\rangle}^{\operatorname{wl}_b=2})^{\operatorname{wl}=n}=
\begin{cases}
\frac{n}{2}-1 & \textrm{for } n\equiv 0 \mod 4\\
\frac{n-1}{2} & \textrm{for } n\equiv 1,3 \mod 4\\
\frac{n}{2} & \textrm{for } n\equiv 2 \mod 4
\end{cases}
\end{align*}
\begin{table}[h]
\centering \caption{Homogeneous basis of $\mathbb{L}_{\langle a,b\rangle}^{\operatorname{wl}_b=2}$}\label{table01}
\begin{center}
\begin{tabular}{@{\hspace{1mm}}c@{\hspace{1mm}}|@{\hspace{1mm}}c
@{\hspace{1mm}}|@{\hspace{1mm}}c@{\hspace{1mm}}}
$\operatorname{wl}$ & basis & $\dim$ \\
\hline
$1$ & $0$  & $0$  \\
$2$ & $[b,b]$ & $1$ \\
$3$ & $[b,[a,b]]$ & $1$  \\
$4$ & $[b,[a,[a,b]]]$ & $1$ \\
$5$ & $[b,[a,[a,[a,b]]]$, \quad $[[a,b],[a,[a,b]]]$ & $2$ \\
$6$ & $[b,[a,[a,[a,[a,b]]]]]$,\quad $[[a,b],[a,[a,[a,b]]]]$, & $3$ \\& $[[a,[a,b]],[a,[a,b]]]$ & \\
$7$ & $[b,[a,[a,[a,[a,[a,b]]]]]]$,\quad $[[a,b],[a,[a,[a,[a,b]]]]]$,& $3$ \\&  $[[a,[a,b]],[a,[a,[a,b]]]]$ & \\
$8$ & $[b,[a,[a,[a,[a,[a,[a,b]]]]]]]$, \quad$[[a,b],[a,[a,[a,[a,[a,b]]]]]]$, & $3$ \\& $[[a,[a,b]],[a,[a,[a,[a,[a,b]]]]]]$ & \\
$9$ & $[b,[\overbrace{a,[a,[\ldots,[a}^{7},b] \ldots ]$,\quad $[[a,b],[\overbrace{a,[a,[\ldots,[a}^{6},b] \ldots ]$, & $4$ \\&
$[[a,[a,b]],[\overbrace{a,[a,[\ldots,[a}^{5},b]\ldots]$,\quad $[[a,[a,[a,b]]],[a,[a,[a,[a,b]]]]]$&   \\
$10$ & $[b,[\overbrace{a,[a,[\ldots,[a}^{8},b] \ldots ]$,\quad $[[a,b],[\overbrace{a,[a,[\ldots,[a}^{7},b] \ldots ]$, & $5$ \\&
$[[a,[a,b]],[\overbrace{a,[a,[\ldots,[a}^{6},b]\ldots]$,\quad $[[a,[a,[a,b]]],[\overbrace{a,[a,[\ldots,[a}^{5},b]\ldots]$, &  \\&
$[[a,[a,[a,[a,b]]]],[a,[a,[a,[a,b]]]]]$ &  \\
$11$ & $[1,10]$, $[2,9]$, $[3,8]$,$[4,7]$,$[5,6]$ & $5$\\
$12$ & $[1,11]$, $[2,10]$, $[3,9]$,$[4,8]$,$[5,7]$ & $5$\\
$\vdots$ & $\vdots$ \\
$n$ & $[1,n-2]$, $[2,n-2]$, \ldots
\end{tabular}
\end{center}
\end{table}
\end{prop}
\begin{prf}In order to prove this we run Algorithm \ref{alg01} with the following choices for the $x_i$ for $i\geq 1$---the complements $W_i$ will always be chosen to be generated by remaining constructed basis elements. Set $x_1:=a$ (and $W_1:=\langle b\rangle$). This results in $Z_1=\langle b\rangle \oplus \langle [a,b]\rangle \oplus \langle [a,a]\rangle$. Now set $x_2:=[a,a]$ (and $W_2:=\langle b,[a,b]\rangle$); clearly, $\deg x_2$ is even. Since $[[a,a],b]=2\cdot [a,[a,b]]$, it follows that
\begin{align*}
Z_3&=W_2\oplus [[a,a],W_2]\oplus [[a,a],[[a,a],W_2]]\oplus \ldots
\\&= \bigoplus_{j\geq 0} \langle \ad^j([a,a])(b)\rangle \oplus \bigoplus_{j\geq 0} \langle \ad^j([a,a])([a,b])\rangle
\\&= \bigoplus_{j\geq 0} \langle [\overbrace{a,[a,[\ldots,[a}^{\operatorname{wl}_a=2j \textrm{ is even}},b] \ldots ]\rangle \oplus \bigoplus_{j\geq 0} \langle [\overbrace{a,[a,[\ldots,[a}^{\operatorname{wl}_a=2j+1 \textrm{ is odd}},b] \ldots ] \rangle
\\&= \bigoplus_{j\geq 0} \langle [\overbrace{a,[a,[\ldots,[a}^{j},b] \ldots ] \rangle
\\&=\bigoplus_{j\geq 0} \langle \ad^j(a)(b)\rangle
\end{align*}
With $x_3:=b$ (of odd degree) we obtain that
\begin{align*}
Z_4&=\langle [b,b]\rangle \oplus \bigoplus_{j>0} \langle \ad^j(a)(b)\rangle \oplus \bigoplus_{j>0} [b[\overbrace{a,[\ldots,[a}^{j},b] \ldots ] \rangle
\\&= \underbrace{\bigoplus_{j> 0} \langle \ad^j(a)(b)\rangle}_{(*)} \oplus \underbrace{\bigoplus_{j\geq 0} \langle [b, \ad^j(a)(b)]\rangle}_{(**)}
\end{align*}
We continue this algorithm with
\begin{align}\label{eqnxi}
x_{3+i}:=[\overbrace{a,[a,[\ldots,[a}^{i},b] \ldots ]=\ad^i(a)(b)
\end{align}
for successive $i\geq 1$.
The terms in $(*)$ have word length $1$ in $b$, the ones in $(**)$ already have word length $2$ in $b$. As we are only interested in a homogeneous basis of $\mathbb{L}_{\langle a,b\rangle}^{\operatorname{wl}_b=2}$ and since word length is additive under brackets, we observe that the chosen $x_i$ for $i\geq 4$ are sufficient to inductively construct such a basis---splitting off any other constructed basis element necessarily leaves a basis for word length $2$ in $b$ untouched, as brackets formed with it in the algorithm will have word length at lest $3$ in $b$.

Consequently, it remains to elicit the effect of successively splitting off the $x_i$ for $i\geq 4$ on $\mathbb{L}_{\langle a,b\rangle}^{\operatorname{wl}_b\leq 2}$---these are the elements which under bracketing still may affect $\mathbb{L}_{\langle a,b\rangle}^{\operatorname{wl}_b=2}$---and hence on $\mathbb{L}_{\langle a,b\rangle}^{\operatorname{wl}_b=2}$ itself, in particular. That is, from now on we shall only construct the spaces $Z'_i$ for $i\geq 5$ consisting of those newly constructed basis elements from $Z_i$ of word length at most $2$ in $b$---i.e.~they will denote the intersection $Z_i':=Z_i \cap \mathbb{L}_{\langle a,b\rangle}^{\operatorname{wl}_b\leq 2}$. We do the same with the chosen canonical complement of the $x_i$, i.e.~$W_i':=W_i \cap \mathbb{L}_{\langle a,b\rangle}^{\operatorname{wl}_b\leq 2}$.

We make the following observations concerning the algorithm:
\begin{itemize}
\item
By construction, the partial sum $\bigoplus_{j> i-4} \langle \ad^j(a)(b)\rangle$ from $(*)$ is a vector subspace of $Z_i$ for $i\geq 4$. This guarantees that we indeed can successively split of the $x_i$ as defined in \eqref{eqnxi}. Moreover, in the end these are the only basis elements of word length $1$ in $b$, i.e.~
\begin{align*}
\mathbb{L}_{\langle a,b\rangle}^{\operatorname{wl}_b=1}=\langle \ad^j(a)(b)\rangle_{j\geq 0}
\end{align*}
\item
By construction, the space $(**)$ is a vector subspace of any $Z_i$ for $i\geq 4$. In other words, the elements $[b, \ad^j(a)(b)]$ for $j>0$ are basis elements of $\mathbb{L}_{\langle a,b\rangle}^{\operatorname{wl}_b=2}$.
\item Any non-trivial Lie product of $x_i$ for $i\geq 4$ with $(**)$ has word length $3$ in $b$ and is not in any $Z'_i$ for $i\geq 4$ nor in $\mathbb{L}_{\langle a,b\rangle}^{\operatorname{wl}_b=2}$.
\item Note that for $i\geq 4$ we have that $\deg x_i \equiv i \mod 2$, since $\deg a$, $\deg b$ are odd. Hence, adapting Algorithm \ref{alg01}, we obtain for $i\geq 4$ that
\begin{align*}
Z'_{i+1}=\mathbb{L}_{\langle a,b\rangle}^{\operatorname{wl}_b\leq 2} \cap
\begin{cases}
W_i'\oplus \langle [x_i,W_i'] \rangle \oplus \langle [x_i,x_i] \rangle & \textrm{if $i$ is odd}\\
W_i' \oplus \langle [x_i,W_i'] \rangle & \textrm{if $i$ is even}
\end{cases}
\end{align*}
with basis elements $[b, \ad^j(a)(b)]\in W_i'$ for all $i,j>0$.
\item This amounts to the fact that passing from $Z_i'$ to $Z_{i+1}'$ we split off one generator $x_i$, and we add new basis elements to $Z_{i+1}'$ of the form
\begin{align*}
\begin{cases}
[x_i,x_j]_{j>i}, [x_i,x_i] & \textrm{if $i$ is odd}\\
[x_i,x_j]_{j>i} & \textrm{if $i$ is even}
\end{cases}
\end{align*}
\item Hence ``summing'' over all degrees $i\geq 1$, and collecting the split-off elements $x_i$, we obtain a basis of $\mathbb{L}_{\langle a,b\rangle}^{\operatorname{wl}_b\leq 2}$ of the form
\begin{align*}
[x_k,x_l]_{k\geq l\geq 3}, (x_k)_{k\geq 3}, [a], [a,a]
\end{align*}
where, by abuse of notation, we implicitly ignore the $[x_i,x_i]=0$ for $i$ even. With the shorthand notation introduced afore the Proposition, namely $[i,j]=[x_{i+2},x_{j+2}]$,
this yields a basis of $(\mathbb{L}_{\langle a,b\rangle}^{\operatorname{wl}_b=2})^{\operatorname{wl}=i}$ which is exactly as asserted.
\end{itemize}
The dimension formula follows trivially. This finishes the proof.
\end{prf}

We again draw on $\deg a$, $\deg b$ being odd. The subsequent corollary is crucial in the construction of the counter-examples, as it will guarantee that the infinitely many elements we newly add to $\mathbb{L}_{\langle a,b\rangle}$, and which possess a non-trivial differential a priori, really do produce new relations and not new spherical homology, i.e., in particular, no new homotopy.
\begin{cor}\label{cor02}
For each odd natural number $j\geq 1$ there exists an element $y_{j} \in (\mathbb{L}_{\langle a,b\rangle}^{\operatorname{wl}_b=2})^{\deg=i_j}$ of odd degree
\begin{align*}
i_j=2\cdot \deg b + j\cdot \deg a
\end{align*}
and word length $\operatorname{wl}(y_j)=j+2$ with the property that
\begin{align*}
y_j \not\in I(y_{k})_{k<j, k \textrm{ odd}}
\end{align*}
where $I(y_{k})_{k<j, k \textrm{ odd}}$ denotes the ideal generated by all the elements $y_1,y_3,y_5,\ldots,\linebreak[4] y_{j-2}$
within $\mathbb{L}_{\langle a,b\rangle}^{\operatorname{wl}_b=2}$
\end{cor}
\begin{prf}
First, note that an element in $\mathbb{L}_{\langle a,b\rangle}^{\operatorname{wl}_b=2}$ has word length $j+2$ if and only if it has degree $i_j$. This allows us to pass from degree to word length and back in the following.

We construct the $y_j$ inductively, starting with, $y_1:=[b,[a,b]]\neq 0$, $y_3,y_5,\ldots$. Suppose we have constructed all suitable $y_1, y_3,y_5,\ldots, y_{j-2}$ (over odd indices). The ideal $I(y_{k})_{k<j, k \textrm{ odd}}$ generated by these elements in $\mathbb{L}_{\langle a,b\rangle}^{\operatorname{wl}_b=2}$ is the vector subspace
\begin{align*}
I(y_{k})_{k<j, k \textrm{ odd}}=\langle \ad^{l_1}(a)(y_1),\ldots , \ad^{l_{j-1}}(a)(y_{j-1}) \rangle_{l_1,\ldots, l_{j-1}\geq 0}
\end{align*}
resulting from successively multiplying all $y_k$ (for odd $k<j$) with all elements of word length $0$ in $b$---recall that, inductively, it holds that  $\operatorname{wl}_{b}(y_k)=2$.

Indeed, for this we made the trivial observation (see Lemma \ref{lemma01}) that $(\mathbb{L}_{\langle a,b\rangle}^{\operatorname{wl}_b=0})=\langle a, [a,a]\rangle$ (basically since $[a,[a,a]]=0$). Using the Jacobi identity we deduce
\begin{align*}
[[a,a],x]=2[a,[a,x]]
\end{align*}
(since $\deg a$ is odd and irrespective of the parity of the degree of $x$). This (together with (anti-)commutativity) implies the asserted form of $I(y_{k})_{k<j, k \textrm{ odd}}$.

In word length $n\geq \operatorname{wl}(y_{j-2})=j$ we specify this to
\begin{align*}
(I(y_{k})_{k<j, k \textrm{ odd}})^{\operatorname{wl}=n}&=\langle \ad^{n-\operatorname{wl}(y_1)}(a)(y_1),\ldots , \ad^{n-\operatorname{wl}(y_{j-2})}(a)(y_{j-2}) \rangle
\intertext{respectively, in degree $n\geq \deg y_{j-2}=2\cdot\deg b + (j-2)\cdot\deg a$, this is}
(I(y_{k})_{k<j, k \textrm{ odd}})^{\deg=n}&=\langle \ad^{\tfrac{n-\deg(y_1)}{\deg a}}(a)(y_1),\ldots , \ad^{\tfrac{n-\deg(y_{j-2})}{\deg a}}(a)(y_{j-2}) \rangle
\end{align*}

We trivially deduce the respective vector space dimensions
\begin{align}\label{eqn02}
(I(y_{k})_{k<j, k \textrm{ odd}})^{\operatorname{wl}=n},\dim (I(y_{k})_{k<j, k \textrm{ odd}})^{\deg=n} = (j-1)/2
\end{align}
(as the $k,j$ are odd).
In particular, the same holds true in degree $n=i_j=2\cdot\deg b + j\cdot\deg a$ in which we aim to find the new element $y_j$ with $\deg y_j=i_j$.

\bigskip

Due to Proposition \ref{prop01} we have that
\begin{align}\label{eqn03}
\dim (\mathbb{L}_{\langle a,b\rangle}^{\operatorname{wl}_b=2})^{\operatorname{wl}=n}=\frac{n-1}{2}
\end{align}
for odd $n$, which is clearly strictly monotonously increasing in $n$, and which yields that $\dim (\mathbb{L}_{\langle a,b\rangle}^{\operatorname{wl}_b=2})^{\operatorname{wl}=j+2}=\frac{j+1}{2}$. Combining Properties \eqref{eqn02} and \eqref{eqn03}, it follows that
\begin{align*}
\dim (\mathbb{L}_{\langle a,b\rangle}^{\operatorname{wl}_b=2})^{\deg=i_{j}}
=\frac{j+1}{2} >
\frac{j-1}{2} =\dim (I(y_{k})_{k<j, k \textrm{ odd}})^{\deg=i_j}
\end{align*}
Consequently, we may pick a non-trivial element
\begin{align*}
y_j \in (\mathbb{L}_{\langle a,b\rangle}^{\operatorname{wl}_b=2})^{\deg=i_{j}} \setminus (I(y_{k})_{k<j, k \textrm{ odd}})^{\deg=i_j}
\end{align*}
(necessarily of word length $j+2$). This finishes the induction step and the proof.
\end{prf}

Hence we may provide the following construction pattern drawing on Corollary \ref{cor02} and the notation therein.
\begin{constr}\label{constr01}
We construct the minimal free connected cochain Lie algebra $(\mathbb{L}_{\langle a, b, x_j\rangle_{j\geq 1, j \textrm{ odd}}}, \delta)$. Choose two generators $a$, $b$ of positive odd degree and form the free graded Lie algebra $(\mathbb{L}_{\langle a,b\rangle},0)$ with $\delta|_{\langle a,b \rangle}=0$. We extend this algebra by an infinite sequence of additional generators $(x_j)_{j\geq 1}$ of strictly monotonously increasing even degrees (the $j=1,3,\ldots$ are chosen to be odd natural numbers)
\begin{align*}
\deg x_j=i_j+1=2\cdot \deg b + j\cdot \deg a+1
\end{align*}
In Corollary \ref{cor02} we constructed a specific sequence of elements
\begin{align*}
y_j\in (\mathbb{L}_{\langle a,b\rangle}^{\operatorname{wl}_b=2})^{\deg=i_j}
\end{align*}
This allows us to specify  well-defined differentials
\begin{align*}
\delta x_j=y_j
\end{align*}
and completes the construction of $(\mathbb{L}_{\langle a,b,x_j\rangle_{j\geq 1, j \textrm{ odd}}},\delta)$. We use spatial realisation (see Section \ref{subsec03})
\begin{align*}
X=|\langle C^*(\mathbb{L}_{\langle a, b, x_j\rangle_{j\geq 1, j \textrm{ odd}}}, \delta)\rangle|
\end{align*}
to construct the topological space $X$. As discussed, we may extend $X$ for example via products or one-point unions.
\end{constr}
\begin{rem}\label{rem02}
Beside choosing arbitrarily positive odd degrees of $a$ and $b$ and extending by products or one-point unions we have further freedom in constructing similar Lie algebras. Indeed, we may pick \emph{any} infinite subsequence of the $x_i$, and the line of argument in the subsequent proof of Theorem \ref{theoA} will just apply analogously. Moreover, we then also have much more freedom to choose the differentials of the $x_i$.
\end{rem}

We are now in the position to give the
\begin{proof}[\textsc{Proof of Theorem \ref{theoA}}]
Construction \ref{constr01} clearly yields a minimal free connected cochain Lie algebra $(\mathbb{L}_{\langle a, b, x_j\rangle_{j\geq 1, j \textrm{ odd}}}, \delta)$ of finite type and a CW-complex $X$ as its spatial realisation. It remains to verify the following properties of $X$.

\begin{enumerate}
\item Even-degree cohomology is finite dimensional, but odd-dimensional cohomology is infinte dimensional, and the space is \linebreak[4]$(\min\{\deg a,\deg b\})$-connected.
\item Spherical cohomology is finite dimensional.
\item The space is coformal.
\end{enumerate}

\bigskip

\noindent \textbf{ad (1).} We use the fact that $(\mathbb{L}_{\langle a, b, x_j\rangle_{j\geq 1, j \textrm{ odd}}}, \delta)$ is a minimal free Lie model of $X$ (see Section \ref{subsec03}) by construction. Applying Properties \eqref{eqn06} we hence obtain that
\begin{align*}
(\langle a, b, x_j\rangle_{j\geq 1, j \textrm{ odd}})_{i-1}=H_{i-1}(\langle a, b, x_j\rangle_{j\geq 1, j \textrm{ odd}},\delta_{\langle a, b, x_j\rangle_{j\geq 1, j \textrm{ odd}}})\cong H_{i}(X)
\end{align*}
for $i\geq 2$. Since $\deg a$, $\deg b$ are odd, and $\deg x_j$ is always even for the chosen odd $j\geq 1$, we deduce that
\begin{align*}
\dim H_\even(X;\qq)&=\dim \langle a,b\rangle=2<\infty
\intertext{and}
\dim H_\odd(X;\qq)&=\dim \langle x_j\rangle_{j \geq 1, j \textrm{ odd}}=\infty
\end{align*}
Moreover, the smallest non-trivial degree of $X$ up to suspension is given by the degree of the smallest generator in even degrees, i.e.~the space is
$(\min\{\deg a,\deg b\})$-connected. In particular, it is always simply-connected. Since we can pick $\deg a$ and $\deg b$ arbitrarily large, the space $X$ is arbitrarily highly connected.

\bigskip

\noindent\textbf{ad (2).} We have seen that the space has \emph{finite dimensional even-degree rational cohomology}, but \emph{infinite dimensional rational cohomology in odd-degrees}. Once we have seen that its \emph{rational spherical homology is finite dimensional}, $X$ constitutes a counter-example to Conjecture \ref{conj01}. As for the latter, it remains to see that none of the chosen $x_j$ constitutes spherical cohomology---in fact, we use that they are separated by degree.

\bigskip

In Remark \ref{rem01} we observed that $v\in \langle a, b, x_j\rangle_{j\geq 1, j \textrm{ odd}}$ is spherical if and only if
\begin{align*}
\delta v \in \delta([\mathbb{L}_{\langle a, b, x_j\rangle_{j\geq 1, j \textrm{ odd}}},\mathbb{L}_{\langle a, b, x_j\rangle_{j\geq 1, j \textrm{ odd}}}]) \In \mathbb{L}_{\langle a, b, x_j\rangle_{j\geq 1, j \textrm{ odd}}}
\end{align*}

In Corollary \ref{cor02} we proved that any considered
\begin{align}\label{eqn07}
\delta x_j=y_j \not\in I(y_{k})_{k<j}
\end{align}
does not lie in the ideal generated by the $y_k$ in $\mathbb{L}_{\langle a,b\rangle}$ for odd $k<l$. Since $\delta$ is a differential, and since $\delta(x_j)\in \mathbb{L}_{\langle a,b\rangle}$, by word length considerations in the $x_j$---the differential $\delta$ decreases word length in the $x_j$ by at most one---it follows that
\begin{align*}
y_j=\delta x_j \in \delta([\mathbb{L}_{\langle a, b, x_j\rangle_{j\geq 1, j \textrm{ odd}}},\mathbb{L}_{\langle a, b, x_j\rangle_{j\geq 1, j \textrm{ odd}}}])
\end{align*}
if and only if $y_j \in I(y_{k})_{k<j}\In \mathbb{L}_{\langle a,b\rangle}$. In view of \eqref{eqn07}, this proves that none of the $x_j$ constitutes spherical cohomology, i.e.~
\begin{align*}
\operatorname{Spher}_*(X)=s\langle a,b\rangle
\end{align*}
of finite dimension $2$.

\bigskip

\noindent\textbf{ad (3).} We have proved that $X$ indeed is a counter-example to the Omnibus conjecture. The additional property of being coformal now follows from the existence of the quasi-isomorphism of differential graded Lie algebras
\begin{align*}
(\mathbb{L}_{\langle a,b,x_j\rangle_{j\geq 1, j \textrm{ odd}}},\delta) \to (H(\mathbb{L}_{\langle a,b,x_j\rangle_{j\geq 1, j \textrm{ odd}}},\delta),0)
\end{align*}
mapping $a\mapsto [a]$, $b\mapsto [b]$ and all $x_j\mapsto 0$. This morphism induces an isomorphism on cohomology basically due to the very same arguments as we applied to compute spherical homology above.
\end{proof}

\begin{rem}
In order to compute the dimensions $\dim (\mathbb{L}_{\langle a,b\rangle}^{\operatorname{wl}_b=i})^j$ instead of using Algorithm \ref{alg01} one may as well dualise the Lie algebra to a Sullivan algebra. For commutative differential graded algebras \textsc{Sage 9.0} provides efficient tools for computing cohomology. This then results in the following scheme of dimensions given in Table \ref{table02} (where word length runs on the vertical axis, word length in $b$ on the horizontal axis).
\begin{table}[h]
\centering \caption{Dimension of $(\mathbb{L}_{\langle a,b\rangle}^{\operatorname{wl}_b=i})^j$}\label{table02}
\begin{center}
\begin{tabular}{@{\hspace{2mm}}c@{\hspace{2mm}}|@{\hspace{2mm}}c
@{\hspace{2mm}} @{\hspace{2mm}}c@{\hspace{2mm}} @{\hspace{2mm}}c@{\hspace{2mm}}@{\hspace{2mm}}c@{\hspace{2mm}}@{\hspace{2mm}}c@{\hspace{2mm}}@{\hspace{2mm}}c@{\hspace{2mm}}
@{\hspace{2mm}}c@{\hspace{2mm}}@{\hspace{2mm}}c@{\hspace{2mm}}@{\hspace{2mm}}c@{\hspace{2mm}}|@{\hspace{2mm}}c@{\hspace{2mm}}}
$j\backslash i$ & $0$ & $1$ & $2$ &$3$ & $4$ & $5$ & $6$ & $7$ & $8$ & $\Sigma=\dim (\mathbb{L}_{\langle a,b\rangle})^j$\\
\hline
$1$ & $1$ & $1$ & $0$ &$0$ & $0$ & $0$ & $0$ & $0$ & $0$ & $2$\\
$2$ & $1$ & $1$ & $1$ &$0$ & $0$ & $0$ & $0$ & $0$ & $0$ & $3$\\
$3$ & $0$ & $1$ & $1$ &$0$ & $0$ & $0$ & $0$ & $0$ & $0$ & $2$\\
$4$ & $0$ & $1$ & $1$ &$1$ & $0$ & $0$ & $0$ & $0$ & $0$ & $3$\\
$5$ & $0$ & $1$ & $2$ &$2$ & $1$ & $0$ & $0$ & $0$ & $0$ & $6$\\
$6$ & $0$ & $1$ & $3$ &$3$ & $3$ & $1$ & $0$ & $0$ & $0$ & $11$\\
$7$ & $0$ & $1$ & $3$ &$5$ & $5$ & $3$ & $1$ & $0$ & $0$ & $18$\\
$8$ & $0$ & $1$ & $3$ &$7$ & $8$ & $7$ & $3$ & $1$ & $0$ & $30$\\
$9$ & $0$ & $1$ & $4$ &$9$ & $14$ & $14$ & $9$ & $4$ & $1$ & $56$\\
\end{tabular}
\end{center}
\end{table}
The obvious ``duality'' in the table clearly just results from exchanging $a$ and $b$. Note that the total dimensions $\dim (\mathbb{L}_{\langle a,b\rangle})^j$ grow exponentially in $j$ due to hyperbolicity. We are already familiar with column $i=2$ due to Table \ref{table01}.
\end{rem}


\section{``Dual'' results}\label{sec04}

Already in the introduction, motivated by the duality of free connected cochain Lie algebras and Sullivan algebras, we brought up the following question ``dual'' to Conjecture \ref{conj01}.
\begin{ques}\label{ques01}
Let $X$ be a simply-connected space with rational homology of finite type. If $\dim \pi_\odd(X)\otimes \qq<\infty$ and $\dim \operatorname{Spher}^*(X)<\infty$. Does then $\dim \pi_*(X)\otimes \qq<\infty$ hold true?
\end{ques}
As already indicated in the introduction, let us quickly discuss known results on the Omnibus Conjecture and their relations to this question, before, in Example \ref{ex01}, we shall give a simple counter-example to it.
\begin{itemize}
\item
Baues' theorem (see \cite{Bau77}, cf.~\cite[Example, p.~328]{FHT01}) shows that the Omnibus Conjecture is true whenever $H_\even(X;\qq)=\qq$. Indeed, in this case a minimal free Lie model of $X$ is concentrated in even degrees; consequently, for degree reasons, its differential vanishes. It follows that the assumption of finite spherical homology implies that the Lie model is finitely generated, i.e.~that the rational homology of $X$ is finite dimensional (see Section \ref{subsec03} and Equations \eqref{eqn06}).
\item Passing from the free Lie model to the Sullivan model we may easily mimic this result: Question \ref{ques01} holds true if $\pi_\odd(X)\otimes \qq=0$. In this case, a minimal Sullivan model is concentrated in even degrees, and finite dimensional spherical cohomology implies finite dimensional rational homotopy groups $\dim \pi_*(X)\otimes \qq<\infty$.
\item The Omnibus Conjecture holds true if $X$ is formal. In this case $(H^*(X;\qq),0)$ encodes the rational type of $X$, so every algebra generator of $H^*(X;\qq)$ in odd degrees is necessarily a spherical cohomology class. Hence, odd-degree rational cohomology is infinite dimensional if and only if there are infinitely many odd-degree spherical cohomology classes.
\item We dualise this situation and assume that $X$ is additionally coformal. Then, however, we obtain the negative answer provided by Example \ref{ex01} below. This is basically due to the fact that a finitely generated free Lie algebra may be infinite dimensional irrespective of the degrees it is generated in, in contrast to a Sullivan algebra, which is finite dimensional if it is finitely generated in odd degrees.
\end{itemize}

Let us now answer Question \ref{ques01} in the negative with
\begin{ex}\label{ex01}
Consider the minimal Sullivan algebra in normal form (see Section \ref{subsec02}) given by
\begin{align*}
V=\langle a,b,n_1,n_2,n_3, \ldots\rangle
\end{align*}
with $\deg a=3$, $\deg b=2$, $\deg n_1=4$, $\deg n_2=6$, \ldots, $\deg n_i=2i+2$ and with well-defined differentials specified by
\begin{align*}
\dif a=\dif b=0,\quad \dif n_1=a\cdot b, \quad \dif n_i=a\cdot n_{i-1} \textrm{ for } i\geq 2
\end{align*}
Again, we may apply spatial realisation in order to obtain the CW complex $X=|\langle(\Lambda V,\dif)\rangle|$ (see Section \ref{subsec03}). It follows directly that
\begin{align*}
\dim \pi_\odd(X)\otimes \qq&=\dim V^\odd=\dim \langle a \rangle=1<\infty\\
\dim \pi_*(X)\otimes \qq&>\dim \pi_\even(\Lambda V,\dif)=\dim V^\even= \dim \langle b, n_i\rangle_{i\geq 1}=\infty\\
\dim \operatorname{Spher}^*(X)&=\dim \langle a,b\rangle=2<\infty
\end{align*}
Note that in this case ``dual'' to the Omnibus Conjecture it is trivial to see (just by word length and the fact that $n_i$ is the only generator with differential involving $n_{i-1}$ for $i\geq 2$) that the $n_i$ do \emph{not} constitute spherical cohomology (see Section \ref{subsec02}).
Consequently, finite dimensional spherical cohomology together with finite dimensional odd homotopy does not imply finite dimensional homotopy groups, which easily answers Question \ref{ques01} in the negative.

Note that again this one example obviously gives rise to an infinite number of arbitrarily highly connected examples $X$ by preserving the stucture but starting with $1\leq \deg a=\deg b+1$ an arbitrarily odd natural number.

Observe further that the example clearly is coformal (due to quadratic differentials only), yet not formal due to the existence of the non-trivial cohomology class $0\neq [6bn_1n_2-2n_1^3-3b^2n_2^2]$.
\end{ex}

\begin{ex}
It is even simpler to construct an example of a space $X$ with finite dimensional spherical cohomology, $\dim \pi_\even(X)\otimes \qq<\infty$, but infinite dimensional $\pi_\odd(X)\otimes \qq$. For this just consider the free Lie algebra $(\mathbb{L}_{\langle a,b\rangle},0)$ on two even-degree generators realised by the CW complex $X$. In this case spherical cohomology up to duality is given by $\langle a, b\rangle$, $\pi_\even(X)\otimes \qq=0$, $\dim \pi_\odd(X)\otimes \qq=\infty$, as $(\mathbb{L}_{\langle a,b\rangle},0)$ is concentrated in even degrees (see Section \ref{subsec03} and Equations \eqref{eqn06}).

Geometrically speaking, this Lie algebra is realised by the one-point union of two odd-dimensional spheres $\s^{\deg a+1}\vee \s^{\deg b+1}$. This is a rationally hyperbolic space (i.e.~it is simply-connected, has finite dimensional rational cohomology, but infinite dimensional rational homotopy), and the Lie algebra is necessarily infinte dimensional.

An example ``dual'' to this would be a space $X$ with finite dimensional spherical homology, $\dim H^\odd(X)<\infty$, but $\dim H^*(X;\qq)=\infty$. Trivially, we may just take $X=\cc\pp^\infty$ for this.
\end{ex}


\section{The elliptic Omnibus Conjecture}\label{sec06}

The author is grateful to Aniceto Murillo and Yves F\'elix for independently speculating about the alternative and stronger version of the Omnibus Conjecture which we documented in Conjecture \ref{conj02} and which we named the ``elliptic Omnibus Conjecture''. Clearly, our counter-examples do not apply to this stronger version of the conjecture, as their rational homotopy is infinite.

One may be tempted to ``dualise'' this conjecture as well leading to
\begin{ques}\label{ques02}
Is every simply-connected topological space $X$ of finite type with finite dimensional rational cohomology and finite dimensional \emph{odd-degree} rational homotopy rationally elliptic?
\end{ques}
Note that $\cc\pp^\infty$ is a non-elliptic space with $\dim \pi_*(\cc\pp^\infty)\otimes \qq=1<\infty$, and vanishing \emph{odd-degree} cohomology. Likewise, $\s^3\vee \s^3$ is a hyperbolic space with $\dim H^*(\s^3\vee \s^3)=3$. Its free minimal Lie model is $(\mathbb{L}_{\langle a,b\rangle},0)$ with $\deg a=\deg b=2$ whence its rational homotopy is concentrated in odd degrees. This should justify the parities chosen in the question.

\bigskip

Concluding the article, it is our goal to illustrate in the following how intimately the (elliptic) Omnibus Conjecture is related to yet another prominent problem in rational homotopy theory, namely Anick's conjecture.

Although or maybe since elliptic spaces are far from being generic (even in the class of spaces satisfying only one of the two finiteness conditions) it is a very interesting question to which extent an arbitrary space may ``differ'' from such an elliptic one. One question in this direction is whether it is in some sense possible to ``embed'' a space with one of the finiteness conditions into an elliptic space.

Anick conjectured that any simply-connected
finite CW-complex $S$ can be realised as the $k$-skeleton of some elliptic complex as
long as $k > \dim S$, or, equivalently, any simply-connected finite Postnikov piece
$S$ can be realized as the base of a fibration $F\hto{} E\to S$ where $E$ is elliptic and $F$ is
$k$-connected with $k > \dim S$.

In \cite{FJM04} this famous problem in rational homotopy theory (see \cite[Chapter 39, Problem 3, p.~516]{FHT01}) was formulated as
\begin{conj}[Anick]\label{conj05}
Any simply connected finite CW complex $S$ can be approximated
arbitrarily closely on the right by an elliptic space. That is, for each natural
number $n$ there is an elliptic space $E_n$ and an $n$-equivalence $S\to E_n$.

Equivalently, any simply connected finite Postnikov piece $S$ can be approximated
arbitrarily closely on the left by an elliptic space. That is, for each
natural number $n$ there is an elliptic space $E_n$ and an $n$-equivalence $E_n \to S$.
\end{conj}
We recall from \cite[Conjecture A, p.~3]{FJM04} that by passing to minimal Sullivan models Anick's conjecture (in the formulation for Postnikov towers) clearly can be rephrased in terms of Sullivan models.
\begin{conj}\label{conj03}
Suppose given a minimal Sullivan algebra $(\Lambda W,\dif)$ with $\dim W<\infty$, $W^1=0$, and fix $n\in \nn$. Then there is an elliptic relative Sullivan algebra $(\Lambda W\otimes \Lambda V,\dif)$ which is even minimal as a non-relative Sullivan algebra and which satisfies that $V=V^{\geq n}$.
\end{conj}

In the same article the following conjecture is presented and then shown to imply Anick's conjecture (see \cite[Proposition 5, p.~6]{FJM04}). We took the freedom to present this conjecture by F\'elix--Jessup--Murillo actually in a slightly weaker form (replacing an arbitrary KS-basis by one induced by ordinary upper degree) and to label it.
\begin{conj}[nilpotent Omnibus Conjecture]\label{conj04}
Let $(\Lambda V,\dif)$ be a minimal Sullivan algebra, $V^1=0$, in which both $V^\even$ and
$H^\even(\Lambda V,\dif)$ are finite dimensional. Then, given a homogeneous basis $(v_i)_{i\geq 0}$ of $V$,
each $v_i$ is nilpotent in the cohomology algebra $H(\Lambda V^{\geq \deg v_i},\bar \dif)$.
\end{conj}
Let us quickly recall the reasoning providing that Anick's conjecture (in the form of Conjecture \ref{conj03}) follows. Indeed, given some minimal Sullivan algebra $(\Lambda W,\dif)$ as in \ref{conj03}, by introducing possibly infinitely many \emph{odd degree} generators with suitable differentials in the minimal algebra from some degree $N\geq 3$ upwards, i.e.~by forming the minimal algebra $(\Lambda(W\oplus Z^\odd),\dif)$, one may guarantee that the outcome has finite dimensional even degree cohomology. The nilpotent Omnibus Conjecture applies and yields that every element $w_i\in W^\even$ is nilpotent in $H(\Lambda W^{\geq i},\dif)$. (Clearly, any such element represents a cohomology class in the quotient $\Lambda W^{\geq i}=\Lambda W/\Lambda W^{<i}$ explaining the slight abuse of notation in the conjecture.) Thus, since $\dim W^\even<\infty$ the same actually holds true when adding only finitely many odd degree generators to $W$ in the first place, i.e.~a finite dimensional homogeneous subspace $\tilde Z\In Z$. Using the subsequent proposition derived from \cite[Theorem 4, p.~5]{FJM04} it follows that $(\Lambda W\otimes \Lambda \tilde Z,\dif)$ is elliptic.

\bigskip

Let $(\Lambda V,\dif)$ be a minimal Sullivan algebra with $\dim V<\infty$, $V^1=0$, and denote by $(v_i)_i$ a homogeneous basis of $V$.
\begin{prop}\label{propfin}
The following properties are equivalent.
\begin{itemize}
\item $\dim H(\Lambda V,\dif)<\infty$,
\item Each $[v_i]$ is nilpotent in the cohomology of the fibre of the rational fibration
\begin{align*}
(\Lambda V^{<\deg v_i},\dif) \hto{} (\Lambda V,\dif) \to (\Lambda V^{\geq \deg v_i}, \bar \dif)
\end{align*}
\end{itemize}
\end{prop}

\bigskip

We extend this picture of the interdependencies between the different conjectures as follows.
\begin{prop}\label{propomn}
For each respective space we have the following implications between the conjectures.
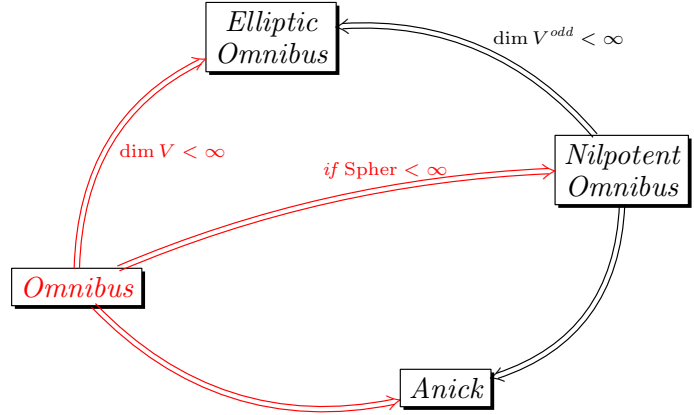
\begin{align*}
\xymatrix{
& *+[F-,]\txt{Elliptic\\Omnibus } +<5pt> && \\
&&& *+[F-,]\txt{Nilpotent\\Omnibus}  \ar@/_2pc/@2{->}[llu]_(0.43){\txt{\tiny \quad ~ $\dim V^\odd<\infty$}} \ar@/^2pc/@2{->}[ldd] \\
*+[F-,]\txt{\textcolor{red}{Omnibus}} \ar@/^2pc/@2{->}@[red][ruu]_(0.43){\txt{\textcolor{red}{\tiny $\dim V<\infty$}}}  \ar@/_2pc/@2{->}@[red][rrd]  \ar@/^1pc/@2{->}@[red][rrru]^(0.6){\txt{\textcolor{red}{\tiny if $\operatorname{Spher}<\infty$}}} &&& \\
&& *+[F-,]\txt{Anick} &
}
\end{align*}
Whenever $\dim V<\infty$, all the Omnibus Conjectures agree.
\end{prop}
A posteriori speaking, in view of our counter-examples to the Omnibus Conjecture the implications drawn in red are a mere triviality when speaking about the conjectures \emph{in general}. However, clearly, also an a priori reasoning \emph{for each separate space} is possible as presented in the proof.

Note that the Omnibus Conjecture and the nilpotent Omnibus Conjecture speak about different classes of spaces, and that there are spaces the second conjecture applies to in contrast to the first one. In order to make sense of the implication between them we restrict to the class of spaces dealt with in the Omnibus Conjecture by additionally imposing that spherical cohomology be finite dimensional.
\begin{proof}[\textsc{Proof of Proposition \ref{propomn}}]
The Omnibus Conjecture clearly implies its special case given by the elliptic one.

As for the Omnibus Conjecture implying the nilpotent one, we suppose the situation of the latter. As announced, however, in order to be able to apply the Omnibus Conjecture we additionally assume that the spherical cohomology of $(\Lambda V,\dif)$ is finite dimensional. Hence we derive that $\dim H(\Lambda V,\dif)<\infty$ due to the Omnibus Conjecture. By Proposition \ref{propfin} we derive the nilpotent one.

Due to Proposition \ref{propfin} the nilpotent and the elliptic conjecture are equivalent whenever $\dim V<\infty$.

We recalled that the nilpotent conjecture implies Anick's conjecture afore the statement of the proposition.

The last implication follows by transitivity. Indeed, the reasoning providing that the nilpotent Omnibus Conjecture implies Anick's Conjecture adds high degree generators to the model. It thereby only creates relations and no new spherical cohomology, i.e.~the construction is performed in such a way that the newly constructed model is in normal form (see Section \ref{subsec02}), and no elements of the underlying vector space lie in $\ker \dif$. Hence, indeed, the Omnibus Conjecture applies and provides the identical conclusions as the nilpotent one.
\end{proof}

Clearly both Anick's conjecture and the elliptic Omnibus Conjecture deal with elliptic spaces. Moreover, note that the elliptic Omnibus Conjecture actually forms an obstruction to performing the construction (recalled above) which derives Anick's conjecture from the nilpotent Omnibus Conjecture. That is, it grants that there is always a non-trivial cohomology class in suitably high even degrees (unless cohomology is finite dimensional) and that the process of successively adding odd-degree generators with suitable differentials can be sustained. It then would remain to prove that it suffices to add only finitely many odd degree generators in order to guarantee finite dimensional cohomology.

Hence, in view of this and with the difference between Conjectures \ref{conj02} and \ref{conj04} mainly being so to say an ``odd-degree-finiteness'' property, it is probably fair to say that the elliptic Omnibus Conjecture is intimately related to and can be considered a (maybe only slightly weaker) interesting variation of Anick's Conjecture.



\begin{thebibliography}{1}

\bibitem{Bau77}
H.~J. Baues.
\newblock Rationale {H}omotopietypen.
\newblock {\em Manuscripta Math.}, 20(2):119--131, 1977.

\bibitem{FHT01}
Y.~F{\'e}lix, S.~Halperin, and J.-C. Thomas.
\newblock {\em Rational homotopy theory}, volume 205 of {\em Graduate Texts in
  Mathematics}.
\newblock Springer-Verlag, New York, 2001.

\bibitem{FHT15}
Y.~F{\'e}lix, S.~Halperin, and J.-C. Thomas.
\newblock {\em Rational homotopy theory. {II}}.
\newblock World Scientific Publishing Co. Pte. Ltd., Hackensack, NJ, 2015.

\bibitem{FJM04}
Y.~Felix, B.~Jessup, and A.~Murillo-Mas.
\newblock Anick's conjecture for spaces with decomposable {P}ostnikov
  invariants.
\newblock {\em Math. Proc. Cambridge Philos. Soc.}, 137(3):559--570, 2004.

\bibitem{Nei10}
J.~Neisendorfer.
\newblock {\em Algebraic methods in unstable homotopy theory}, volume~12 of
  {\em New Mathematical Monographs}.
\newblock Cambridge University Press, Cambridge, 2010.

\end{thebibliography}

\def\cprime{$'$}


\vfill

\begin{center}
\noindent
\begin{minipage}{\linewidth}
\small \noindent \textsc
{Manuel Amann} \\
\textsc{Institut f\"ur Mathematik}\\
\textsc{Differentialgeometrie}\\
\textsc{Universit\"at Augsburg}\\
\textsc{Universit\"atsstra\ss{}e 14 }\\
\textsc{86159 Augsburg}\\
\textsc{Germany}\\
[1ex]
\footnotesize
\textsf{manuel.amann@math.uni-augsburg.de}\\
\textsf{www.uni-augsburg.de/de/fakultaet/mntf/math/prof/diff/team/dr-habil-manuel-amann/}
\end{minipage}
\end{center}

\end{document}